\documentclass [a4paper, 12t]{article}
\usepackage{amsmath}
\usepackage{amsfonts}
\usepackage{amssymb}
\usepackage{amsthm}
\usepackage{makeidx}
\usepackage{graphicx}
\usepackage{pb-diagram}
\usepackage{epstopdf}
\usepackage{fancyhdr}
\usepackage{enumitem}
\usepackage[all]{xy}

\newtheorem{teo}{Theorem}

\newtheorem{lem}{Lemma}
\newtheorem{cor}{Corollary}
\theoremstyle{definition}

\newtheorem{oss}{Observation}

\newtheorem{con*}{Conjecture}

\newtheorem{rem} {Remark}
\begin{document}
\setlength{\parskip}{1ex plus 0.5ex minus 0.2ex}
\begin{center}

\textbf{ A CHARACTERIZATION OF SEPARABLE CONJUGATE BANACH SPACES}\\
\end{center}
\begin{center}
\textsl{By STEFANO ROSSI}
\end{center}
$$$$
$$$$
\textbf{Abstract}\quad The following elegant characterization of dual Banach spaces will be proved: a \emph{separable} Banach space $\mathfrak{X}$ is a dual space if and only exists a norm closed subspace $\mathfrak{M}\subset\mathfrak{X}^*$ with (Dixmier) characteristic equal to $1$, whose functionals are \emph{norm-attaining}.
\tableofcontents
\begin{section} {Introduction}
Among all Banach spaces, dual spaces have some additional properties, which assure a more far-reaching treatment of their structure. Probably the most useful property is represented by weak*-compactness of the unit ball of such spaces. Giving necessary and sufficient conditions for a Banach space $\mathfrak{X}$ to be a dual space is a long standing problem; it goes back essentially to J. Dixmier \cite{DixmierDuke}. Many theorems concerning this intriguing  area are known; these results usually require the compactness property  for $\mathfrak{X}_1$, the unit ball of $\mathfrak{X}$, with respect to suitable locally convex topology. Elegant as theoretical results, these theorems, however, do not provide any practical tool to solve the question in the applications, because compactness property is very difficult to be checked in infinite-dimensional contexts.\\
In spite of this, reflexive Banach spaces (which are indeed dual spaces with \emph{unique} predual) are well understood thanks to a celebrate thorem by R. C. James \cite{James1}, which states that a Banach space $\mathfrak{X}$ is reflexive if and only if every continuous linear functional attains its norm on $\mathfrak{X}_1$.\\
The result quoted above can be regarded as a straightforward application of the following powerful characterization of weak compactness, which is due as well to James:
\begin{teo} [James, \cite{James2}]
\label{huge}
Let $\mathfrak{Y}$ be a \emph{complete} locally convex space and $\mathcal{C}\subset\mathfrak{Y}$ a bounded, $\sigma(\mathfrak{Y},\mathfrak{Y}^*)$ closed subset The following are equivalent:
\begin{enumerate}
\item $\mathcal{C}$ is a weakly compact subset.
\item Every $\varphi\in \mathfrak{Y}^*$ attains its sup on $C$,i.e. for every $\varphi\in\mathfrak{Y}^*$ there exists $c\in \mathcal{C}$ such that $|\varphi(c)|=\sup_{y\in\mathfrak{C}}|\varphi(y)|$.
\end{enumerate}
\end{teo}
The proof of the theorem announced in the abstract essentially depends on James ' theorem \ref{huge}, where completeness assumption unfortunately cannot be removed, as James himself has shown in \cite{Jamesb}.
\end{section}
\begin {section} {Notations and preliminaries}
If $\mathfrak{X}$ is a Banach space, $\mathfrak{X}_1$ will stand for its unit ball, that is $$\mathfrak{X}_1=\{x\in\mathfrak{X}: \|x\|\leq 1\}$$ The dual (conjugate) space of $\mathfrak{X}$ will be denoted by $\mathfrak{X}^*$. A Banach space $\mathfrak{X}$ is a dual space if there exists a Banach space $\mathfrak{Y}$ such that $\mathfrak{Y}^*\cong\mathfrak{X}$ (isometric isomorphism); in this case $\mathfrak{Y}$ is said to be a \emph{predual}.\\
Accordig to Dixmier, we define the characteristic of a subspace $\mathfrak{M}\subset\mathfrak{X}^*$ as the the real number given by
$$c(\mathfrak{M})\doteq\inf_{x\neq 0}\left\{\sup_{\varphi\in\mathfrak{M}_1}\frac{|\varphi(x)|}{\|x\|} \right\}$$
 When $c(\mathfrak{M})=1$, we say that $\mathfrak{M}$ is a \emph{determinant} subspace of the conjugate space; a standard application of the bipolar theorem shows that $\mathfrak{M}$ is determinant iff $\mathfrak{M}_1$ is weakly* dense in $\mathfrak{X}^*_1$.\\
Finally, we recall that a continuous linear functional $\varphi\in\mathfrak{X}^*_1$ is said to be \emph{norm-attaining} if there exist $x\in\mathfrak{X}_1$ such that $|\varphi(x)=\|\varphi\|$. A subspace $\mathfrak{M}\subset \mathfrak{X}^*$ is norm-attainig if each $\varphi\in \mathfrak{M}$ is a norm-attaining linear functional.
\end {section}
\begin{section} {A general result}
This section is devoted to present a general result for the existence of preduals. In what follows, we will need some tools from  the general theory of locally convex space, in particular some basic facts belonging to the theory of dual pair. For an extensive treatment see, for istance, \cite{Wil}. For completeness' sake, here we recall that, given a dual pair $(\mathfrak{X},\mathfrak{Y})$, the Mackey topology on $\mathfrak{X}$ is the locally convex topology described by the base of seminorms $p_K$ given by
$$p_K(x)=sup_{y\in K}|\langle x, y\rangle|$$
where $K$ is a $\sigma(\mathfrak{Y},\mathfrak{X})$-compact, (circled) convex subset of $\mathfrak{Y}$. The Mackey topology is often denoted by $\tau(\mathfrak{X},\mathfrak{Y})$. Analogously, one can define the Mackey topology $\tau(\mathfrak{y},\mathfrak{X})$ on $\mathfrak{Y}$.\\

 We begin our analysis by observing that any conjugate space $\mathfrak{X}^*$ is a complete locally convex space with the Mackey topology $\tau(\mathfrak{X}^*,\mathfrak{X})$. This can be proved, for instance, with the aid Grothendieck completeness theorem or as an easy application of the lemma \ref{compactconv} from section \emph{4} of the present paper. A second observation needed is that any predual of a Banach space $\mathfrak{M}$ should be sought as a suitable subspace of the conjugate space $\mathfrak{X}^*$. More precisely, one has the following:
 \begin{lem}
 Let $\mathfrak{X}$ and $\mathfrak{M}$ be Banach spaces.
 If the isometric isomorphism $\mathfrak{M}^*\cong \mathfrak{X}$ holds, then there is an isometric injection $i:\mathfrak{M}\rightarrow\mathfrak{X}^*$ such that $i(\mathfrak{M})$ is a determinant subspace, whose functionals are norm-attaining.
 \end{lem}
\begin{proof}
Let $\Phi:\mathfrak{X}\rightarrow\mathfrak{M}^*$ be an isometric isomorphism. Taking the (Banach) adjoint of $\Phi$, we get an isometric isomorphism $\Phi^*:\mathfrak{M}^{**}\rightarrow \mathfrak{X}^*$. Let $i:\mathfrak{M}\rightarrow\mathfrak{X}^*$ the linear isometry given by the composition $\Phi^*\circ j$, where $j$ is the canonical injection of $\mathfrak{M}$ into its second conjugate space $\mathfrak{M}^{**}$. It only remains to prove that $i(\mathfrak{M})$ has the required properties.\\
 Let $x\in \mathfrak{M}$, then we have $$\|x\|=\|\Phi(x)\|=\sup_{m\in \mathfrak{M}_1} |\langle\Phi(x), m\rangle|=\sup_{m\in\mathfrak{M}_1} |\langle i(m),x\rangle|$$ The last equality says that $i(\mathfrak{M})$ is determinant.\\In order to prove that $i(\mathfrak{M})$ is norm attaining, observe that, given any $m\in\mathfrak{M}$, one has $\|i(m)\|=\|m\|=|\varphi(m)|$, for some $\varphi\in \mathfrak{M}^*_1$ by virtue of the Hanh-Banach theorem. Now, if $x$ be is the unique element of $\mathfrak{X}$ such that $\Phi(x)=m$, we have $\|i(m)\|=\langle i(m), x\rangle$. This concludes the proof, since $\|x\|=\|\Phi(x)\|=\|m\|\leq 1$.
\end{proof}
Under suitable additional requests the previous lemma can be reverted to get existence of a predual, as indicated in the following result (due to Dixmier \cite{DixmierDuke}), which can be regarded as a partial convers of Alaoglu theorem. For the completeness' sake, we provide a different prove based on the Krein-Smulian theorem.
\begin{teo}
\label{compact}
Let $\mathfrak{M}\subset\mathfrak{X}^*$ be a norm closed determinant subspace, such that $\mathfrak{X}_1$ is $\sigma(\mathfrak{X},\mathfrak{M})$- compact. Then $\mathfrak{M}^*\cong\mathfrak{X}$, through the isometric isomorphism $\Phi:\mathfrak{X}\rightarrow\mathfrak{M}^*$ given by $\langle\Phi(x), \varphi\rangle=\varphi(x)$ for all $\varphi\in \mathfrak{M}^*$ and for all $x\in\mathfrak{X}$.
\end{teo}
\begin{proof}
The linear map $\Phi:\mathfrak{X}\rightarrow\mathfrak{M}^*$ described in the statement is surely isometric, since $ \mathfrak{M}$ is determinant. Let $\mathfrak{N}\subset\mathfrak{M}^*$ be its range. $\mathfrak{N}$ is weakly* dense in $\mathfrak{M}^*$, since it is a total subspace, \emph{i.e.} $\mathfrak{N}^{\perp}=0$ (in the duality $(\mathfrak{M}^*,\mathfrak{M})$).  To get the conclusion, it is enough to observe that $\mathfrak{N}$ is weakly* closed by virtue of the Krein-Smulian theorem (here is the point werw completeness of $ \mathfrak{M}$ is needed). In fact $\mathfrak{N}_1=\Phi(X_1)$ is weakly* compact as the continuous image of a compact set ($\Phi$ is apparently a $\sigma(\mathfrak{X},\mathfrak{M})-\sigma(\mathfrak{M}^*,\mathfrak{M})$ continuous map).
\end{proof}
Actually the theorem remains true under the milder assumption that $\mathfrak{M}$ is just a total subspace, see \cite{Kaijser}, anyway
Theorem \ref{compact} is enough to state and prove the following general characterization of dual Banach space:
\begin{teo}
\label{mainteo}
Let $\mathfrak{X}$ be a Banach space. The following conditions are equivalent:
\begin{enumerate}
\item $\mathfrak{X}$ is a dual space.
\item There exists a norm closed subspace $\mathfrak{M}\subset \mathfrak{X}^*$ which is determinant, norm-attaining and such that $\mathfrak{X}$ is complete with respect to the Mackey topology $\tau(\mathfrak{X},\mathfrak{M})$.
\end{enumerate}
Moreover, each subspace as in 2. is canonically\footnote{ $\mathfrak{M}^*\cong\mathfrak{X}$ through the isometric isomorphism $\Phi:  \mathfrak{X}\rightarrow\mathfrak{M}^*$ given by the restriction on $\mathfrak{M}$ of the canonical injection $j:\mathfrak{X}\rightarrow\mathfrak{X}^{**}$, i.e. $\langle\Phi(x),\varphi\rangle=\varphi(x)$ for each $\varphi\in\mathfrak{M}$.} a predual of $\mathfrak{X}$ .
\end{teo}
\begin{proof}
The implication $\emph{1}.\Rightarrow \emph{2}.$ has been already proved. Conversely, if $\emph{2.}$ holds, we can equip $\mathfrak{X}$ with the Mackey topology $\tau(\mathfrak{X},\mathfrak{M})$. By the Mackey-Arens theorem, the Mackey dual of $\mathfrak{X}$ is $\mathfrak{M}$. The unit ball $X_1\subset X$ is clearly Mackey-bounded since it is even norm bounded; moreover $X_1$ is $\sigma(\mathfrak{X},\mathfrak{M})$-closed. In fact, thanking to the fact that $\mathfrak{M}$ is determinant, we can rewrite $X_1$ as $$\left\{x\in X:|\varphi(x)|\leq 1\quad\textrm{for all}\, \varphi\in \mathfrak{M}_1\right\}$$
so it is a $\sigma(\mathfrak{X},\mathfrak{M})$- closed subset as intersection of $\sigma(\mathfrak{X},\mathfrak{M})$-closed sets. By the assumption of completeness, James'theorem applies, so we get $\sigma(\mathfrak{X},\mathfrak{M})$-compactness of $X_1$. Hence $\mathfrak{M}$ is canonically a predual of $\mathfrak{X}$.
\end{proof}
\begin{rem}
Unfortunately completeness hypothesi is not omissible, as the following counterexample shows. Let us consider the Banach space $l_1[0,1]$ of (real) function $f$ defined on $[0,1]$ such that $\sum_{x\in [0,1]}|f(x)|<\infty$, endowed with the corresponding $\|\cdot\|_1$- norm. Its conjugate space is clearly represented by $B[0,1]$, the space of all bounded function over $[0,1]$ with the sup norm, acting by $\langle g, f\rangle=\sum_{x\in [0,1]}g(x)f(x)$ for all $g\in B[0,1]$, $f\in l_1[0,1]$. The norm closed subspace of continuous $C[0,1]\subset B[0,1]$ is determinant and norm-attaining. However it is not a predual, since $C[0,1]^*\cong\mathcal{M}([0,1])$ (the space of all finite Borel measures on [0,1]), while the inclusion\footnote{One think $f\in l_1[0,1]$ as the Borel measure $\sum_{x\in [0,1]}f(x)\delta_x$.} $l_1[0,1]\subset \mathcal{M}([0,1])$ is proper. Actually the trouble is that $l_1[0,1]$ is not complete under the Mackey topology $\tau(l_1[0,1], C[0,1])$, as pointed out in \cite{Josè}.
\end{rem}
Although elegant, the previous theorem is not completely satisfactory, because to handle completeness notions in the general context of locally convex spaces is rather difficult; however the result can be considerably strengthened in the \emph{separable case}, as it is shown in the next section.
\end{section}
\begin{section} {The separable case}
In this section we devote our analysis to \emph{separable} Banach spaces $\mathfrak{X}$, the reason being that separability assumption allows us to prove the following useful completeness result:
\begin{teo}
\label{completeness}
Let $\mathfrak{X}$ be a separable Banach space and $\mathfrak{M}\subset \mathfrak{X}^*$ a determinant subspace. $X$ is a \emph{complete} locally convex space with respect to the Mackey topology $\tau(\mathfrak{X},\mathfrak{M})$.
\end{teo}
The proof of theorem \ref{completeness} relies on the following two lemmas:
\begin{lem}
\label{compactconv}
Let $\mathfrak{X}$ be topological vector space. If $\{f_{\alpha}:\alpha\in I\}\subset X^*$ is a net of continuous linear functionals converging to $f$ uniformly over the compact subset of $\mathfrak{X}$, then $f$ is \emph{sequentially}-continuous.
\end{lem}
\begin{proof}
It is a standard $\frac{\varepsilon}{3}$-argument.
Let $\{x_n:n\in\mathbb{N}\}\subset X$ be a sequence convergin to $x\in X$. The set $K\doteq\{x_n:n\in\mathbb{N}\}\cup\{x\}$ is compact, thus $f_{\alpha}$ converges to $f$ uniformly on $K$. This means that, given $\varepsilon>0$, there exists $\alpha_0\in I$ such that
$$\sup_{y\in K}|f_{\alpha_0}(y)-f(y)|<\frac{\varepsilon}{3}$$
By the triangle inequality we get
$$|f(x_n)-f(x)|\leq \varepsilon$$
for all $n\geq N_{\varepsilon}$, where $N_{\varepsilon}$ is any natural number such that $|f_{\alpha_0}(x_n)-f_{\alpha_0}(x)|\leq\frac{\varepsilon}{3}$ for all $n\geq N_{\varepsilon}$.
\end{proof}
The second lemma deals with the \emph{Mazur} property of $\mathfrak{M}$ with respect the $\sigma(\mathfrak{M},\mathfrak{X})$-topology. More precisely, we have the following statement:
\begin{lem}
\label{Mazur}
Let $\mathfrak{X}$ be a \emph{separable} Banach space. If $\mathfrak{M}\subset{X}^*$ is a determinant subspace, then $\mathfrak{M}$ equipped with the $\sigma(\mathfrak{M},\mathfrak{X})$-topology is a \emph{Mazur} space, that is every sequantially-continuous linear form  $f:\mathfrak{M}\rightarrow\mathbb{C}$ is continuous with respect to the $\sigma(\mathfrak{M},\mathfrak{X})$-topology.
\end{lem}
\begin{proof}
Since $\mathfrak{X}$ is separable, the weak* topology on $\mathfrak{X}_1$ is metrizable as well its restriction to $\mathfrak{M}_1$. Let $f:\mathfrak{M}\rightarrow\mathbb{C}$ be a $\sigma(\mathfrak{M},\mathfrak{X})$ sequentially-continuous linear form. The restriction $f\upharpoonright_{M_1}$ is thus a $\sigma(\mathfrak{M},\mathfrak{X})$- (uniformly) continuous function. Since $\mathfrak{M}$ is determinant, $\mathfrak{M}_1$ is weakly* dense in $\mathfrak{X}^*_1$, so $f\upharpoonright_{M_1}$ extends to a weakly* continuous function $g:\mathfrak{X}^*_1\rightarrow\mathbb{C}$, which is apparently the restricton to $\mathfrak{X}^*_1$ of a linear form $G$ defined on the whole $\mathfrak{X}^*$. By the Krein-Smulian theorem (see, for instance, \cite{Conway}), this linear functional is weakly*-continuous, that is there is a unique $x\in\mathfrak{X}$ such that $G(\varphi)=\varphi(x)$ for all $\varphi\in \mathfrak{X}^*$. This concludes the proof.
\end{proof}
Finally we can prove theorem \ref{completeness}.
\begin{proof}
Let $\{x_{\alpha}:\alpha\in I\}\subset X$ be a $\tau(\mathfrak{X},\mathfrak{M})$- Cauchy net. This means that, given any\footnote{One should note that if $K\subset\mathfrak{M}$ is $\sigma(\mathfrak{M},\mathfrak{X})$-compact set, then $\overline{\textrm{conv}(K)}^{\sigma(\mathfrak{M},\mathfrak{X})}$ is yet $\sigma(\mathfrak{M},\mathfrak{X})$ compact thank to Alaoglu theorem.} $\sigma(\mathfrak{M},\mathfrak{X})$-compact subset $K\subset\mathfrak{M}$ and any $\varepsilon>0$, we have
\begin{equation}
\label{Cauchy}
\sup_{\varphi\in K}|\varphi(x_{\alpha}-x_{\beta})|<\varepsilon
\end{equation}
for all $\alpha, \beta\in I$ such that $\alpha,\beta\succeq\alpha_0$ for some $\alpha_0\in I$. In particular, fixed $\varphi\in \mathfrak{M}$, $\{\varphi(x_{\alpha}):\alpha\in I\}$ is a numerical Cauchy net, so we can define a linear form on $\mathfrak{M}$, say $G$, by $G(\varphi)\doteq\lim_{\alpha} \varphi(x_{\alpha})$. Thanks to inequality \ref{Cauchy}, the limit is uniform over the $\sigma(\mathfrak{M},\mathfrak{X})$-compact subset of $\mathfrak{M}$. By lemma \ref{compactconv}, $G$ is sequentially $\sigma(\mathfrak{M},\mathfrak{X})$-continuous, so it continuous in this topology thanks to the lemma \ref{Mazur}. This means that there is a unique $x\in \mathfrak{X}$ such that $G(\varphi)=\varphi(x)$ for all $\varphi\in\mathfrak{M}$. Taking the limit (with respect to $\beta\in I$) of the inequality \ref{Cauchy}, one easily get $x_{\alpha}\rightarrow x$ in the $\tau(\mathfrak{M,\mathfrak{X}})$-topology. This concludes our proof.
\end{proof}
$$$$
Finally, we can state our main theorem, as announced in the abstract:

\begin{teo}
\label{sep}
Let $\mathfrak{X}$ be a separable Banach space. Every determinant norm closes subspace $\mathfrak{M}\subset\mathfrak{X}^*$, whose functional are norm-attaining, is canonically a predual of $\mathfrak{X}$.
\end{teo}
\begin{proof}
Theorem \ref{mainteo} of the previous secrion applies, since $\mathfrak{X}$ is complete with respect the Mackey topology $\tau(\mathfrak{X},\mathfrak{M})$, by virtue of completeness theorem \ref{completeness}.
\end{proof}
$$$$
\begin{oss}
The assumption that $\mathfrak{M}$ is norm closed cannot be dropped. Let us consider, for instance, the separable Banach space $C[0,1]$ of continuous function on the compact interval $[0,1]$. Let us consider the subspace $\mathfrak{N}\subset C[0,1]^*$, algebraically spanned by the set $\{\delta_x: x\in [0,1]\}$, where $\delta_x$ is the Dirac measure concentrated on $x\in [0,1]$. It is not difficult to check that $\mathfrak{N}$ is a determinant subspace, whose functionals are norm-attaining. Nevertheless $\mathfrak{N}$ cannot be a predual of $C[0,1]$, since the last space has not preduals at all, as it is well known. Observe that $\overline{\mathfrak{N}}$ is not norm-attaining.
\end{oss}
The following  corollaries are almost immediate.
\begin{cor}
Let $\mathfrak{X}$ be a separable Banach space. If $\mathfrak{M}\subset \mathfrak{X}^*$ is a norm closed subspace as in theorem \ref{sep}, then it is a \emph{maximal} subspace with respect the property to be norm closed and norm-attaining.
\end{cor}
\begin{proof}
If $\mathfrak{N}\supseteq\mathfrak{M}$ is norm closed and norm-attaining, then it is (canonically) a predual of $\mathfrak{X}$ as well as $\mathfrak{M}$, hence $\mathfrak{N}=\mathfrak{M}$, since no proper inclusion relationship among preduals are allowed, thanks to the theorems of Alaoglu-Banach and Krein-Smulian.
\end{proof}
In particular we get the following non-trivial result:
\begin{cor}
If $\mathfrak{X}$ is a (separable) Banach space, such that $\mathfrak{X}^*$ is separable, then $j(\mathfrak{X})\subset\mathfrak{X}^{**}$ is a closed maximal norm-attaining subspace.
\end{cor}
As concrete applications of the corollaries stated above, we get, for instance, that $c_0\subset l_{\infty}$ is a closed subspace, which is maximal with respect to the property to be norm-attaining (on $l_1$), since $c_0$ is a (not unique)  predual of $l_1$. \\
The same is true for $\mathcal{K}(\mathfrak{H})\subset\mathcal{B}(\mathfrak{H})$, where $\mathcal{K}(\mathfrak{H})$ is the norm closed ideal of all compact  operators on the \emph{separable} Hilbert space $\mathfrak{H}$, thought as a predual of $\mathcal{S}_1$, the nuclear (trace-class) operators on $\mathfrak{H}$.
\end{section}
\begin{section} {Applications}
To appreciate theorem \ref{sep}, we illustrate two nice applications of it, simplifying the proofs of two classical results.\\
The first is the Riesz-Frechet representation theorem for the dual of a (separable) Hilbert space, the second one is the Riesz representation theorem for the dual spaces of the Lebesgue spaces $L^p$ ($p>1$) in the separable case.\\
\begin{teo} [Riesz-Frechet]
Let $\mathfrak{H}$ be a separable Hilbert space, with inner product $(\cdot,\cdot)$. If $\varphi\in \mathfrak{H}^*$, then there is a unique $y\in\mathfrak{H}$ such that $\varphi(x)=(x,y)$ for each $x\in\mathfrak{H}$.
\end{teo}
\begin{proof}
There is no loss of generality in the assumption that $\mathfrak{H}$ is a \emph{real} Hilber space; in this case the map $\Phi:\mathfrak{H}\rightarrow\mathfrak{H}^*$, given by $\langle\Phi(y), x\rangle= (x,y)$, is \emph{linear} and isometric. Note that $\Phi(\mathfrak{H})\subset\mathfrak{H}^*$ is norm closed, determinant and norm-attaining. This means that $\Phi(\mathfrak{H})$ is canonically a predual of $\mathfrak{H}$, \emph{i.e.} $\Phi$ is surjective.
\end{proof}
The application to duality theory for $L^p$ spaces is more interesting, since the usual proofs requires rather involved argument invoking Clarkson inequality or Radon-Nikodym theorem. We will consider a measure space $(X,\mathfrak{S},\mu)$, with the property that $\mathfrak{S}$ is countably generated: this assures that the corresponding Lebesgue spaces $L^p(X,\mu)$ are separable for all $p\geq 1$ (for the exception of the value $p=\infty$ ). Given $p>1$, we denote by $q>1$ its conjugate exponent, that is the real number $q$ such that $\frac{1}{p}+\frac{1}{q}=1$. With these notations, the classical Riesz theorem can be stated as follows
\begin{teo}
If $p>1$ is a real number, then the map $\Phi:L^q(X,\mu)\rightarrow {L^p(X,\mu)}^*$ given by $\langle \Phi(f), g\rangle =\int_X fg\textrm{d}\mu$ for all $f\in L^q(X,\mu), g\in L^p(X,\mu)$ , is an isometric isomorphism
.\end{teo}
\begin{proof}
Let $\Psi: L^p(X,\mu)\rightarrow {L^q(X,\mu)}^*$ be the linear map give by
$$\langle\Psi(f), g\rangle\doteq \int_X fg\textrm{d}\mu\quad\textrm{for all}\,f\in L^p(X,\mu), g\in L^q(X\mu) $$
By H\"{o}lder inequality one has $\|\psi(f)\|\leq \|f\|_p$ for all $f\in L^p(X,\mu)$. Given $f\in L^p(X,\mu)$,  let us consider the measurable function $g$ given by
$$g(x)=\left\{
\begin{array} {rl}
|f(x)|^{p-2}f(x)\\
0\,\,\textrm{if}\,\,\, f(x)=0
\end{array}
\right.
$$
It is easy to check that $g\in L^q(X,\mu)$ and that $\|g\|_q={\|f\|_p}^{\frac{p}{q}}=\|f\|_p^{p-1}$. Since one has $\langle\Psi(f),g\rangle=\|f\|_p^p$, it follows that $\langle\Psi(f), h\rangle=\|f\|_p$, where $h\doteq\frac{g}{\|g\|_q}\in L^q(X,\mu)_1$. By virtue of the last equality, we conclude that $\|\Psi(f)\|=\|f\|_p$ and the norm of the functional is attained on the function $h$. This means that $\Psi(L^p(X,\mu))$ is a closed subspace of ${L^q(X,\mu)}^*$, whose functionals are norm-attaining. A similar argument as above also shows that $\Psi(L^p(X,\mu))$ is determinant. By theorem \ref{sep}, we argue that $\Psi(L^p(X,\mu))$ is canonically a predual of $L^q(X,\mu)$,\emph{i.e.} the isometric isomorphism ${L^p(X,\mu)}^*\cong L^q(X,\mu)$ holds as in the statement.
\end{proof}
\begin{rem}
Our proof of Riesz representation theorem only depend on H\"{o}lder inequality.
\end{rem}
\end{section}
$$$$
I wish to thank S. Doplicher for some useful discussions during the peparation of the present paper.

$$$$
\textsl{DIP. MAT. CASTELNUOVO, UNIV. DI ROMA LA SAPIENZA, ROME, ITALY}\\
\emph{E-mail address:} \verb"s-rossi@mat.uniroma1.it"
\end{document}